\documentclass[11pt, a4paper]{amsart}
\usepackage{amssymb,amsmath,amsthm}

\usepackage{mathtools} 

\usepackage[numbers]{natbib}
\usepackage{graphicx}

\newtheorem{lemma}{Lemma}[section]

\theoremstyle{remark}

\theoremstyle{definition}

\usepackage{tikz-cd}
\newtheoremstyle{problemstyle}  
        {3pt}                                               
        {3pt}                                               
        {\normalfont}                               
        {}                                                  
        {\bfseries}                 
        {\normalfont\bfseries.}         
        {.5em}                                          
        {}                                                  
\theoremstyle{problemstyle}

\newtheorem{myproblem}{Problem}

\newtheoremstyle{freethm}
 {3pt}
 {3pt}
 {\itshape}
 {}
 {\bfseries}
 {\normalfont\bfseries.}
 {.5em}
 {} 

\theoremstyle{freethm}
\newtheorem{mytheorem}{Theorem}

\title[The Gaussian Monge problem]{Convergence of the vertical gradient flow for the Gaussian Monge problem}

\author[E.~Jansson]{Erik Jansson$^*$}
\email[]{erikjans@chalmers.se}
	
\author[K.~Modin]{Klas Modin$^{*,\dagger}$} 
\email[]{klas.modin@chalmers.se}

\address{$^*$Department of Mathematical Sciences, Chalmers University of Technology \& University of Gothenburg, S--412~96 G\"oteborg, Sweden.} 
    
\thanks{
$^\dagger$Corresponding author.
}

	\subjclass{15A23, 49M99, 	49Q22  }

	\keywords{Gradient flows, Matrix decompositions, Optimal transport.}
\begin{document}

\maketitle              
\begin{abstract}

We investigate a matrix dynamical system related to optimal mass transport in the linear category, namely, the problem of finding an optimal invertible matrix by which two covariance matrices are congruent.
We first review the differential geometric structure of the problem in terms of a principal fiber bundle. 
The dynamical system is a gradient flow restricted to the fibers of the bundle.
We prove global existence of solutions to the flow, with convergence to the polar decomposition of the matrix given as initial data.
The convergence is illustrated in a numerical example.
\end{abstract}

\section{Introduction}
Matrix decompositions are ubiquitous in mathematics. 
An interesting view-point is to think of matrix decompositions as originating from limits of matrix flows.
Famous examples include the Toda flow as described by \citet{Fl1974} and Brockett's diagonalizing flow~\citep{Br1991b}. 
A review of these and many other examples is given by \citet{Mo2017}, where emphasis is put on the connection to optimal mass transport and information geometry.
The approach is based on Riemannian gradient flows, which reveals connections to optimal transport and information geometry. 
Further, it highlights  hidden dynamical features of matrix factorization, for instance  that the $QR$ decomposition, singular value decomposition and polar decomposition  all can be computed by such gradient flows. 

In this paper, we study the polar decomposition of matrices from a geometric perspective. In particular, we show that it is given as the limit of a \emph{vertical gradient flow}. 
The contribution in this paper is a careful analysis of the finite-dimensional vertical gradient flow introduced by \citet{Mo2017}. Our main result is a proof of convergence to a minimizer. 
The convergence of this flow is interesting in itself, as it completes the picture in \citet{Mo2017}, where its ``cousin'' flow -- the \emph{horizontal gradient flow} -- is shown to converge.
But perhaps more interestingly, the result suggests a similar convergence study in the infinite-dimensional case, for the gradient flow corresponding to Brenier's polar decomposition of maps~\citep{Brenier1991}, by using the corresponding infinite-dimensional geometric framework of gradient flows developed by \citet*{BaKaMo2022}.

\medskip

\noindent\textbf{Acknowledgements.}
This work was supported by the Wallenberg AI,
Autonomous Systems and Software Program (WASP) funded by the Knut and Alice
Wallenberg Foundation. This work was also supported by the Swedish Research Council, grant number 2022-03453, and the Knut and Alice Wallenberg Foundation, grant number WAF2019.0201.

\section{The Gaussian Monge problem}
In this section, we briefly describe the theoretical background, state the optimization problem at the center of our study and outline its geometric structure. 
For details, see \citet{Mo2017} and references therein.

Let $\mu_0,\mu_1$ be zero-mean multivariate Gaussian distributions on $\mathbb{R}^n$. 
The \emph{Gaussian Monge problem} consists of finding an invertible linear map $\phi(x) = Ax$, with $A \in \operatorname{GL}(n)$ an invertible $n\times n$ real matrix, that  pushes $\mu_0$ forward to $\mu_1$ while minimizing the functional $J\colon\operatorname{GL}(n) \to \mathbb{R}$,
\begin{equation}
\label{eq:Mongefunc}
J(A) \coloneqq \int_{\mathbb{R}^n} \|x-Ax\|^2 \mu_0,
\end{equation}
where $\|\cdot\|$ denotes the standard Euclidean norm on $\mathbb{R}^n$. 

It is possible to cast the Gaussian Monge problem only in terms of matrices.
Indeed, from the perspective of information geometry (cf.~\citet{Amari2016}), the set of all zero-mean multivariate Gaussian distributions is a \emph{statistical manifold} isomorphic to the set of covariance matrices, i.e., the set $\operatorname{P}(n)$ of symmetric and positive-definite $n\times n$ matrices. 

We identify $\mu_0$ and $\mu_1$ with their covariance matrices, $\Sigma_0$ and $\Sigma_1$, and note that the functional in equation \eqref{eq:Mongefunc} can be rewritten as 
\begin{align}
\label{eq:functional}
J(A) = \operatorname{Tr}(\Sigma_0(I-A)^T(I-A)),
\end{align}
where $\operatorname{Tr}$ denotes the matrix trace and $I$ is the identity matrix. 
The matrix $A$ transforms $\mathbb R^n$ and thereby it transforms $\mu_0$. Under the identification of Gaussian distributions with covariance matrices, $A$ maps $\Sigma_0$ to $A\Sigma_0A^T$. 
Therefore, the Monge problem for multivariate Gaussian distributions can be formulated in terms of covariance matrices. 

\begin{myproblem}[Gaussian Monge problem]
\label{eq:OT}
\begin{equation}
    \begin{split}
    &\min_{A \in \operatorname{GL}(n)} J(A) = \operatorname{Tr}(\Sigma_0(I-A)^T(I-A))\\
    &\text{subject to: } A \in \mathcal{C}(\Sigma_0,\Sigma_1) = \{A \in \operatorname{GL}(n): A\Sigma_0A^T = \Sigma_1 \}
    \end{split}
\end{equation}
\end{myproblem}
To describe the Riemannian geometry of the problem, we begin by equipping $\operatorname{GL}(n)$ with a Riemannian metric $\mathcal G$ defined by
\begin{equation}
\label{eq:metric}
\mathcal{G}_{A}(\dot A,\dot A) = \operatorname{Tr}(\Sigma_0 \dot A^T \dot A) = \operatorname{Tr}(A\Sigma_0 A^T (\dot A A^{-1})^T (\dot A A^{-1})) 	
\end{equation}
where $\dot A \in T_A \operatorname{GL}(n)$. 
The corresponding Riemannian distance function $d\colon  \operatorname{GL}(n) \times \operatorname{GL}(n) \mapsto \mathbb{R}$ is
\begin{equation*}
	d(A_0,A_1) = \operatorname{Tr}\left(\Sigma_0(A_0-A_1)^T (A_0-A_1)\right)^{1/2},
\end{equation*}
so we see that $J(A) = d^2(I,A)$. 

The Gaussian Monge problem is solved by the polar decomposition of matrices, as pointed out by \citet{Brenier1991}.
\begin{mytheorem}
    Let $A \in \operatorname{GL}(n)$ and  $\Sigma_0 \in \operatorname{P}(n)$. 
    Then there exists unique matrices $P \in \operatorname{P}(n)$ and $Q \in \operatorname{O}(n,\Sigma_0) \coloneqq \{Q \in \operatorname{GL}(n): Q\Sigma_0Q^T = \Sigma_0\}$ such that 
    $$
        A = PQ. 
    $$
    Moreover, the matrix $P$ is the unique solution to Problem \ref{eq:OT} when  $\Sigma_1 = \pi(A) = A\Sigma_0 A^T$. 
    \label{th:brenier}
\end{mytheorem}
\subsection{The geometric structure of the Gaussian Monge problem}
\label{seq:geometry}
We now describe the geometric structure of Problem \ref{eq:OT}.  

First, $\operatorname{O}(n,\Sigma_0)$ is a Lie subgroup of $\operatorname{GL}(n)$ with Lie algebra 
$$
\mathfrak{o}(n,\Sigma_0) = \{X \in \mathfrak{gl}(n): X\Sigma_0 + \Sigma_0 X^T=0\}. 
$$
A matrix $A \in \operatorname{GL}(n)$ acts  on $\Sigma \in P(n)$ by the matrix congruence action $A.\Sigma \coloneqq A\Sigma A^T$. 
This transitive action defines a projection  $\pi\colon \operatorname{GL}(n) \to P(n)$ given by 
$$
 \pi(A) = A\Sigma_0 A^T. 
$$
The isotropy subgroup of $\Sigma_0$  with respect to the congruence action is $\operatorname{O}(n,\Sigma_0)$.  
Further, the Riemannian metric $\mathcal{G}$ on $\operatorname{GL}(n)$ is $\operatorname{O}(n,\Sigma_0)$-invariant. 
That is, $\mathcal{G}$ is invariant under right multiplication with matrices in $\operatorname{O}(n,\Sigma_0)$. Indeed, with an arbitrary $Q \in \operatorname{O}(n,\Sigma_0)$,
\begin{align*}
\mathcal{G}_A(\dot A,\dot A) &= \operatorname{Tr}(\Sigma_0 \dot A^T \dot A)  = \operatorname{Tr}(Q \Sigma_0 Q^T \dot A^T \dot A) \\ &= 	\operatorname{Tr}( \Sigma_0  (\dot A Q)^T \dot A Q) = \mathcal{G}_{AQ}(\dot A Q,\dot A Q). 
\end{align*}

For any $A \in \mathcal{C}(\Sigma_0,\Sigma_1)$, the mapping 
\begin{equation}
    \label{eq:iso}
	 \mathcal{C}(\Sigma_0,\Sigma_1) \ni B \mapsto A^{-1} B \in \operatorname{O}(n,\Sigma_0) 
\end{equation}
provides an isomorphism between $\mathcal{C}(\Sigma_0,\Sigma_1)$ and $\operatorname{O}(n,\Sigma_0)$. 
To see this, note that  $\pi(A Q) = \Sigma_1$ for all $Q \in \operatorname{O}(n,\Sigma_0)$, and for any $B \in \mathcal{C}(\Sigma_0,\Sigma_1)$, $ A^{-1}B \in \operatorname{O}(n,\Sigma_0)$  since $A^{-1} B \Sigma_0 B^T A^{-T} = A^{-1}\Sigma_1 A^{-T} = \Sigma_0$.
Furthermore, we have that $\mathcal{C}(\Sigma_0,\Sigma_1) = \pi^{-1}(\Sigma_1)$. 
The sets $\pi^{-1}(\Sigma)$ are fibers of a principal $\operatorname{O}(n,\Sigma_0)$-bundle over $P(n)$,
\begin{equation}
    \label{eq:bundle}
	\begin{tikzcd}
  \operatorname{GL}(n)  \arrow[d,shorten >=1.5ex]
    & \arrow[l,hook',shorten >=1.5ex]  \operatorname{O}(n,\Sigma_0) \\
P(n)& \end{tikzcd}
\end{equation}

The principal bundle structure \eqref{eq:bundle} gives rise to two linear subspaces of each tangent space $T_A \operatorname{GL}(n)$. 
The \emph{vertical distribution} $\operatorname{Ver}_A$ is geometrically the tangent space of the fiber going through $A$, and is determined by computing the kernel of the derivative of $\pi(A)$, denoted $D\pi(A)$. 
The \emph{horizontal distribution} $\operatorname{Hor}_A$ is the distribution transversal to $\operatorname{Ver}_A$ with respect to the metric $\mathcal{G}$. 
Note that the definition of $\operatorname{Ver}_A$ relies solely on the bundle structure. 
Therefore, $\operatorname{Ver}_A$ is independent of the choice of metric.

To compute $\operatorname{Ver}_A$ we first compute $D\pi(A)$. It is, with $\Sigma = \pi(A)$, given by
$$
\dot \Sigma = \dot A A^{-1} A \Sigma_0 A^T + A \Sigma_0 A^T 	(\dot A A^{-1})^T = \dot A  A^{-1}\Sigma + \Sigma (\dot A A^{-1})^T \coloneqq D\pi(A) \cdot \dot A.
$$
Note that $\dot A A^{-1} \in \mathfrak{gl}(n)$. 
A matrix $V$ is thus in the kernel of $D\pi(A)$ if it satisfies 
\begin{align*}
V \Sigma + \Sigma V^T = 0,
\end{align*}
which is equivalent to the condition $V \in \mathfrak{o}(n,\Sigma)$. 
Thus, the tangent space, called the \emph{vertical distribution} at $A$, is given by 
\begin{equation*}
	\operatorname{Ver}_A = \{\dot A = VA \in T_A \operatorname{GL}(n): V \in \mathfrak{o}(n,\Sigma)\}. 
\end{equation*}

To compute $\operatorname{Hor}_A$, we must compute the orthogonal complement of $ \mathfrak{o}(n,\Sigma)$ with respect to $\mathcal{G}$. 
To this end, let $U \in \mathfrak{gl}(n)$ be  arbitrary and note that if $V \in \mathfrak{o}(n,\Sigma)$, then $V\Sigma = X$ is skew-symmetric, so $V = X\Sigma^{-1}$. 

Inserting $U$ and $V$ into Equation \eqref{eq:metric} we obtain 
\begin{align*}
\mathcal{G}_A(VA,UA) = \operatorname{Tr}(A\Sigma_0 A^T V^T U) = \operatorname{Tr}(\Sigma V^T U)	 = \operatorname{Tr}(X^T U).
\end{align*}
The orthogonal complement of the skew-symmetric matrices under the Frobenius inner product (i.e, $\mathcal{G}$ with $\Sigma_0 = I$), is the set of symmetric matrices, denoted by $\operatorname{Sym}(n)$. Thus, the horizontal bundle is 
$$
\operatorname{Hor}_A = \{\dot A = UA \in T_A \operatorname{GL}(n): U \in\operatorname{Sym}(n)\}. 
$$
Note that the horizontal bundle is independent of $\Sigma_0$ 

The \emph{polar cone} $K_\Diamond$ consists of all matrices in $\operatorname{GL}(n)$ connected with the identity by a horizontal geodesic, i.e.,  geodesics $\gamma(t)$ with $\dot\gamma(0) \in \operatorname{Hor}_{\gamma(0)}$ such that $\dot \gamma(t) \in \operatorname{Hor}_{\gamma(t)}$. 

The following lemma (see \citet{Mo2017} for a geometric proof) is vital to understand the importance of $K_\Diamond$. 
\begin{lemma}
\label{lemma}
The restriction of the projection $\pi$ to $K_\Diamond$ is an isomorphism. 
In other words, $K_\Diamond$ is a section of the principal bundle \eqref{eq:bundle}. 
\end{lemma}

In fact, the polar cone \emph{itself} consists of positive definite symmetric matrices (geometrically these elements are not, however, elements of the base space). 
Now, let $A \in \operatorname{GL}(n)$ be arbitrary and set $\pi(A) = \Sigma_1$. 
By Lemma \ref{lemma}, there is a unique matrix $  K_\Diamond \ni P = \left.\pi\right|_{K_\Diamond}^{-1}(\Sigma_1)$. 
It is clear that $P$ and $A$ both are in the same fiber. 
Therefore, $Q = P^{-1}A \in O(n,\Sigma_0)$, by the isomorphism \eqref{eq:iso}.
Thus, as $PQ = A$, we have obtained the polar decomposition. 
 $P$ lies at the intersection of the polar cone and the fiber $\pi^{-1}(A)$.
The uniqueness of the polar decomposition means in geometric terms that there is only one point of intersection between a fiber and the polar cone.

\section{Vertical gradient flow}

Having introduced the geometry of the Gaussian Monge problem, we move to the gradient flow. 
Let $\nabla_\mathcal{G}$ denote the gradient with respect to the Riemannian metric \eqref{eq:metric}. 
Further, we denote by $\operatorname{Proj}_{\operatorname{Ver}}$ the projection onto the vertical distribution. 

\begin{figure}[ht]
    \centering
    \includegraphics[width=0.7\textwidth]{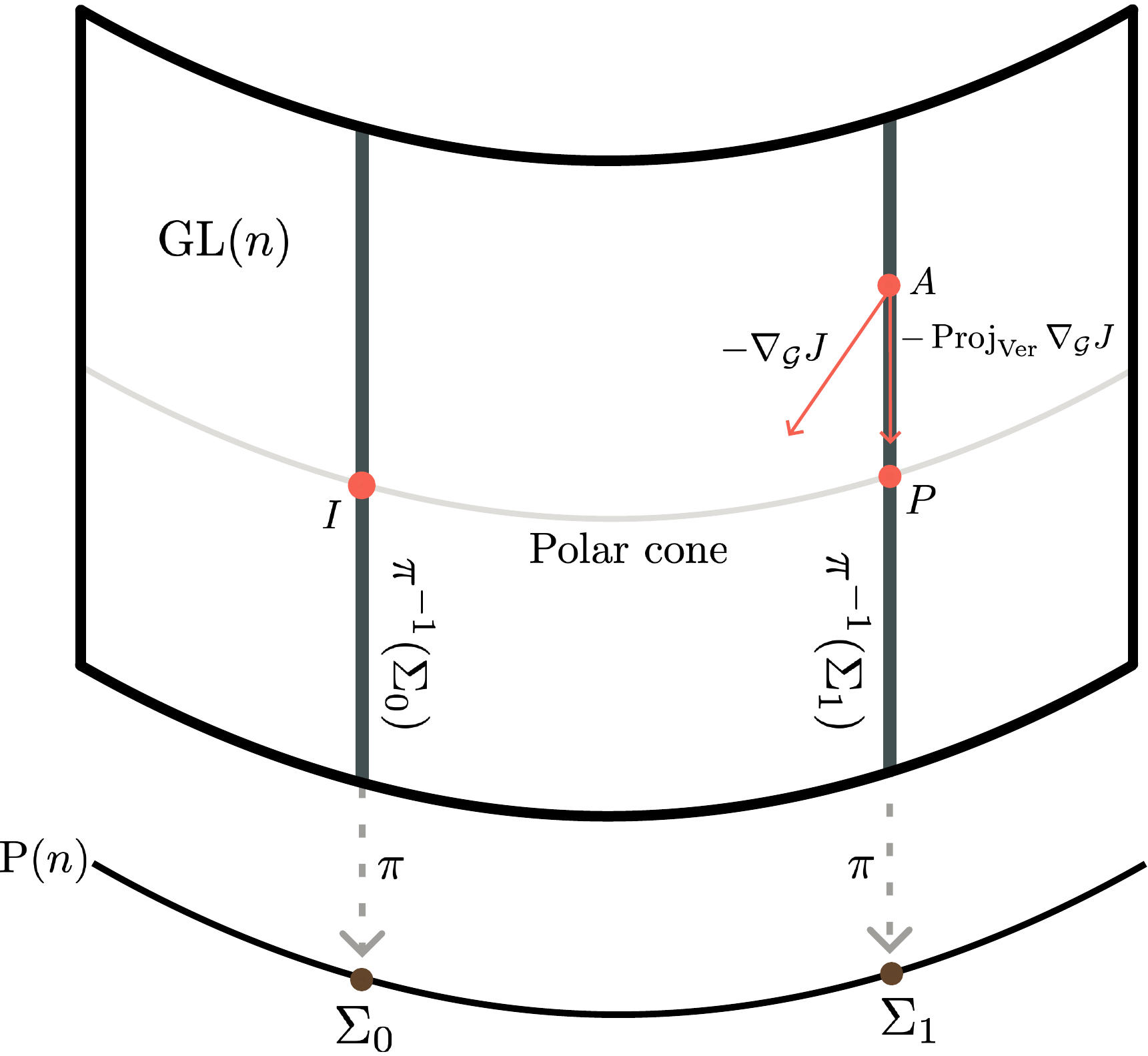}
    \caption{Illustration of the geometry of the problem. Note how the polar cone intersects the fiber at one point, corresponding to the uniqueness of the polar factorization.}
    \label{fig:geom}
\end{figure}

The following vertical gradient flow is suggested in \citet{Mo2017}:
\begin{align}
    \label{eq:gf}
    \dot B = -\operatorname{Proj}_{\operatorname{Ver}} \nabla_\mathcal{G} J(B), ~B(0) = A. 
\end{align}
Note that $\operatorname{Proj}_{\operatorname{Ver}} \nabla_\mathcal{G} J(B)$ is not just the \emph{projection} of the gradient of $J$, but also that it coincides with the gradient with respect to the induced metric on the fiber. 
Therefore, the projected gradient flow is itself a gradient flow on the Riemannian sub-manifold given by the fiber. 
Our goal is to study Equation \eqref{eq:gf}. 
The main result is as follows.
\begin{mytheorem}
    \label{th:main}
    Let $A$ be in the identity component of $\operatorname{GL}(n)$. 
    Then the gradient flow \eqref{eq:gf} has these properties:
    \begin{enumerate}
        \item A global solution exists.
        \item The solution converges as $t \to \infty$ to the matrix $P$ in the polar decomposition of $A$. 
    \end{enumerate}
\end{mytheorem}
For an illustration of the vertical gradient flow as well as the geometry of the Gaussian Monge problem, see Figure~\ref{fig:geom}. 

\subsection{Analysis of the gradient flow}
To prove Theorem \ref{th:main}, we leverage the geometry described in Section \ref{seq:geometry}. 
The first step is to explicitly determine the gradient flow in Equation \eqref{eq:gf}. 
To this end, we take the time derivative of $J(B)$, defined in Equation \eqref{eq:functional}, 
\begin{align*}
    \frac{\mathrm{d}}{\mathrm{d}t} J(B) =2\operatorname{Tr}(\Sigma_0 (B-I)^T \dot B) =  \mathcal{G}_B( 2(B-I), \dot B). 
\end{align*}
Thus, the vertical gradient flow is
\begin{equation}
    \label{eq:gf_lagrange}
    \dot B = -2(B-I) - SB, ~B(0) = A, 
\end{equation}
where $S$ is a Lagrange multiplier forcing $B$ to remain on the fiber.

It is advantageous to rewrite the gradient flow in terms of a right-reduced variable, both from a theoretical and computational perspective. 
Therefore, we introduce the variable $\Omega \coloneqq \dot B B^{-1} $ which, by the definition of the vertical bundle, is in $\mathfrak{o}(n,\Sigma_1)$.
Algebraic manipulations of Equation \eqref{eq:gf_lagrange}, see \citet[Section 2.4.1]{Mo2017} for details, yield the equivalent right--reduced gradient flow, 
\begin{gather}
    \label{eq:gf_RR}
    \dot B = \Omega B, ~B(0) = A, \\
    \Sigma_1 \Omega + \Omega \Sigma_1 = 2\Sigma_1(B^{-1}-B^{-T}). \label{eq:sylvester}
\end{gather}
Finally, note that the right-hand side of Equation \eqref{eq:sylvester} may be expressed without matrix inverses. 
Indeed, it holds that 
\begin{align}
    \label{eq:sav}
    B^{-1} = \Sigma_0 B^T A^{-T}\Sigma_0^{-1}A^{-1}. 
\end{align}
We now prove Theorem \ref{th:main}. 
\begin{proof}
    
        We first prove global existence. 
        We claim that the right--hand side of \eqref{eq:gf_RR} is a locally Lipschitz map.
        To see this, note that it is given by the product of $\Omega$, a function of $B$ as defined by Equation \eqref{eq:sylvester}, and $B$. 
        If we can verify that $B \mapsto \Omega$ is a Lipschitz map, then the claim follows, as the map $B \mapsto B$ is trivially Lipschitz and the product of two Lipschitz maps is a locally Lipschitz map.
        
        The Sylvester equation \eqref{eq:sylvester} is equivalent to the linear system 
        \begin{align}
            \label{eq:linsys}
            \mathcal{S} \operatorname{Vec}(\Omega) = \operatorname{Vec}(2\Sigma_1(B^{-1}-B^{-T})),
        \end{align}
        where $ \mathcal{S} = I \otimes \Sigma_1 + \Sigma_1 \otimes I $ and $\operatorname{Vec}: \mathbb{R}^{n \times n } \to \mathbb{R}^{n^2}$ denotes the vectorization of a matrix. 
        Here $\otimes$ is the Kronecker product of matrices \cite[Chapter 4.3]{HoJo91}. 
        The linear system \eqref{eq:linsys} has a unique solution \citep{BaSt72}, so 
        $\Omega = \operatorname{Vec}^{-1}(\mathcal{S}^{-1}\operatorname{Vec}(2\Sigma_1(B^{-1}-B^{-T})))$. 
       The map $\mathbb{R}^{n^2} \ni x\mapsto \mathcal{S}^{-1}x$ as well as  vectorization and its inverse are Lipschitz maps. Finally,  $B \mapsto 2\Sigma_1(B^{-1}-B^{-T}))$ is Lipschitz since $B^{-1}$ can be computed from $B$  and the fixed matrices $\Sigma_0$ and $A$ using Equation \eqref{eq:sav}. 
       Thus, $B \mapsto \Omega$ is a Lipschitz map. 

       Recall that the gradient flow is on the fiber $\pi^{-1}(\Sigma_1)$ which is isomorphic to the Lie group $\operatorname{O}(n,\Sigma_0)$, on which the metric \eqref{eq:metric} is right-invariant. 
       On a Lie group equipped with a right-invariant metric, locally Lipschitz  guarantees global existence of the flow \citep{BaKaMo2022}.

        To prove convergence, we consider the functional $J(B(t))$ evaluated along the flow as a function of time. 
        A general fact of gradient flows is that the functional is decreasing in time, i.e.,
        \begin{align*}
             \frac{\mathrm{d}}{\mathrm{d}t} J(B) \leq 0. 
        \end{align*}
        We know that $J(B) = d^2(B,I)$ is bounded from below; as $B \in \pi^{-1}(\Sigma_1)$, $J(B)$ must be larger than or equal to $D_{\Sigma_1} \coloneqq d^2(\pi^{-1}(\Sigma_1),\pi^{-1}(I))$, the squared distance between $\pi^{-1}(\Sigma_1)$ and the identity fiber. 
        Therefore, $J(B)$ converges as $t \to \infty$.
        The idea is now to show that $ \frac{\mathrm{d}}{\mathrm{d}t} J(B)$ converges to $0$. 
        By the chain rule one sees that 
        \begin{align}
            \label{eq:dJdt}
           \frac{\mathrm{d}}{\mathrm{d}t} J(B) {=} 2\operatorname{Tr}(\Sigma_0 (B-I)^T \Omega B) {=} 2\operatorname{Tr}(\Sigma_1 \Omega) - 2\operatorname{Tr}(\Sigma_0 \Omega B) {=} {-} 2\operatorname{Tr}(\Sigma_0 \Omega B),
        \end{align}
        where the second-to-last inequality follows by noting that $\pi(B) = \Sigma_1$ and the final since $\Omega\Sigma_1 $ is skew-symmetric and thus has zero trace.
        The second derivative of $J$ is
        \begin{align*}
            \frac{\mathrm{d}^2}{\mathrm{d}t^2} J(B) = - 2\operatorname{Tr}(\Sigma_0 \dot \Omega B)-2\operatorname{Tr}(\Sigma_0  \Omega^2 B),
       \end{align*}
       where $\dot \Omega$ is obtained by solving the Sylvester equation arising from the time derivative of Equation \eqref{eq:sylvester},
    \begin{align*}
        \Sigma_1 \dot \Omega + \dot\Omega \Sigma_1 = 2\Sigma_1  \frac{\mathrm{d}}{\mathrm{d}t} (B^{-1}-B^{-T}).
    \end{align*}
    Note that the Frobenius norm of a solution to a Sylvester equation can be bounded by the Frobenius norm of the right-hand side.
    Indeed, if we let $\operatorname{Sep}(\Sigma_1) \coloneqq \min_{X \in \mathbb{R}^n} \|X\Sigma_1+\Sigma_1 X\|_F/\|X\|_F$, then 
    \begin{align*}
        &\|\dot \Omega\|_F \leq \frac{ \|2\Sigma_1 \frac{\mathrm{d}}{\mathrm{d}t} (B^{-1}-B^{-T})\|_F}{\operatorname{Sep}(\Sigma_1)}, \\
        &\| \Omega\|_F \leq \frac{\| 2\Sigma_1 (B^{-1}-B^{-T})\|_F}{\operatorname{Sep}(\Sigma_1)}.
    \end{align*}
     In our case, $\operatorname{Sep}(\Sigma_1) = 2\lambda_{\min}(\Sigma_1)$ \citep{Varah1979}, where $\lambda_{\min}(\Sigma_1)$ denotes the smallest eigenvalue of $\Sigma_1$. 
     Now, 
     \begin{align*}
        \left |  \frac{\mathrm{d}^2}{\mathrm{d}t^2} J(B) \right|{\leq} 2|\operatorname{Tr}(\Sigma_0 \dot \Omega B)|{+}2|\operatorname{Tr}(\Sigma_0  \Omega^2 B)| {\leq} 2\|\Sigma_0\|_F \left(\|\dot \Omega\|_F{+}\|\Omega\|_F^2\right)\|B\|_F
     \end{align*}
    where the final inequality is due to the Cauchy--Schwarz inequality  and the fact that the Frobenius norm is sub--multiplicative. 
    Recalling Equation \eqref{eq:sav}, it follows by the triangle inequality and the trace-invariance of the Frobenius norm that
    \begin{align*}
        \|\Omega\|_F \leq C_{\Sigma_1,\Sigma_0,A} \|B\|_F,
    \end{align*}
    where $ C_{\Sigma_1,\Sigma_0,A} = 2\|\Sigma_1\|_F \|\Sigma_0\|_F \|\Sigma_0^{-1}\|_F \|A^{-1}\|_F^2/\lambda_{\min}(\Sigma_1)$. 
    Furthermore, Equation \eqref{eq:sav} together with Equation \eqref{eq:gf_RR} gives that 
    \begin{align*}
        \|\dot \Omega\|_F \leq C_{\Sigma_1,\Sigma_0,A}^2\|B\|_F^2,
    \end{align*}
    Therefore,
    \begin{align*}
         \left |  \frac{\mathrm{d}^2}{\mathrm{d}t^2} J(B) \right| \leq 2 C_{\Sigma_1,\Sigma_0,A}^2 \|\Sigma_0\|_F\|B\|_F^3, 
    \end{align*}
    which is finite as $\|B\|_F \leq C\operatorname{Tr}(\Sigma_0 B^T B) = C\operatorname{Tr}(\Sigma_1) <\infty $ for some constant $C>0$. 
    As the second derivative is bounded, Barbalat's lemma \cite[Lemma 8.2]{Kh01} gives that
    $\lim_{t \to \infty} \frac{\mathrm{d}}{\mathrm{d}t} J(B) = 0$. 
    However, the time derivative being zero at some $B \in \pi^{-1}(\Sigma_1)$ means, by Equation \eqref{eq:dJdt}, that 
    \begin{align*}
       0 = \frac{\mathrm{d}}{\mathrm{d}t} J(B) = \operatorname{Tr}(\Sigma_0 \Omega B) = \operatorname{Tr}( \Omega\Sigma_1 B^{-T}). 
    \end{align*}
    where the second inequality is due to the identity $B\Sigma_0 = \Sigma_1 B^{-T}$ and the invariance of the trace under cyclic permutations. 
    As $\Omega \Sigma_1$ by construction is skew-symmetric, $B^{-T}$ must be symmetric, i.e., $B$ is symmetric as well as in the fiber $\pi^{-1}(\Sigma_1)$. 
    From Theorem \ref{th:brenier} it follows that there is only one such matrix, namely $P$. 
    Geometrically, it means that $P$ is the unique point of intersection between the fiber and the polar cone. 
    Thus, the flow $B(t)$ approaches the matrix $P$.  
\end{proof}
\section{A numerical experiment}
In this section, we provide a small numerical example of the gradient flow, illustrating its convergence. For convenience, let $\Sigma_0 = I$. 
To avoid selecting a specific initial matrix $A$, we instead create  matrices of size $2\times 2$ by drawing random matrices with standard normally distributed elements. 
A matrix is discarded if it is not invertible.
The matrix generation process is continued until a set $\{A^i, i =1,\ldots,1000\}$ of $1000$ invertible matrices  is obtained. 

Equation \eqref{eq:gf_RR} is discretized by the \emph{Lie--Euler method} \citep{IsMuKaNoZa2000}, 
\begin{align*}
    B_{k+1}^i = \exp(h \Omega_k^i)B_k^i, \,B_0^i = A^i. 
\end{align*}
Here $\exp$ refers to the matrix exponential, and we set $h = 0.1$ and integrate the flow for $300$ steps.
Further, $\Omega_k$, given as the solution to equation \eqref{eq:sylvester} is computed using the Bartels--Stewart algorithm \citep{BaSt72}. 
We compute the polar decomposition $P^i$  of each $A^i$ from its SVD factorization, as suggested in \cite[Section 3.1]{Hi86}. 

We thus obtain approximated flow curves $B^i_k, k = 1,\ldots, 300, i = 1,\ldots, 1000$. 
 
\begin{figure}[ht!]
    \centering
    \includegraphics[width =\textwidth]{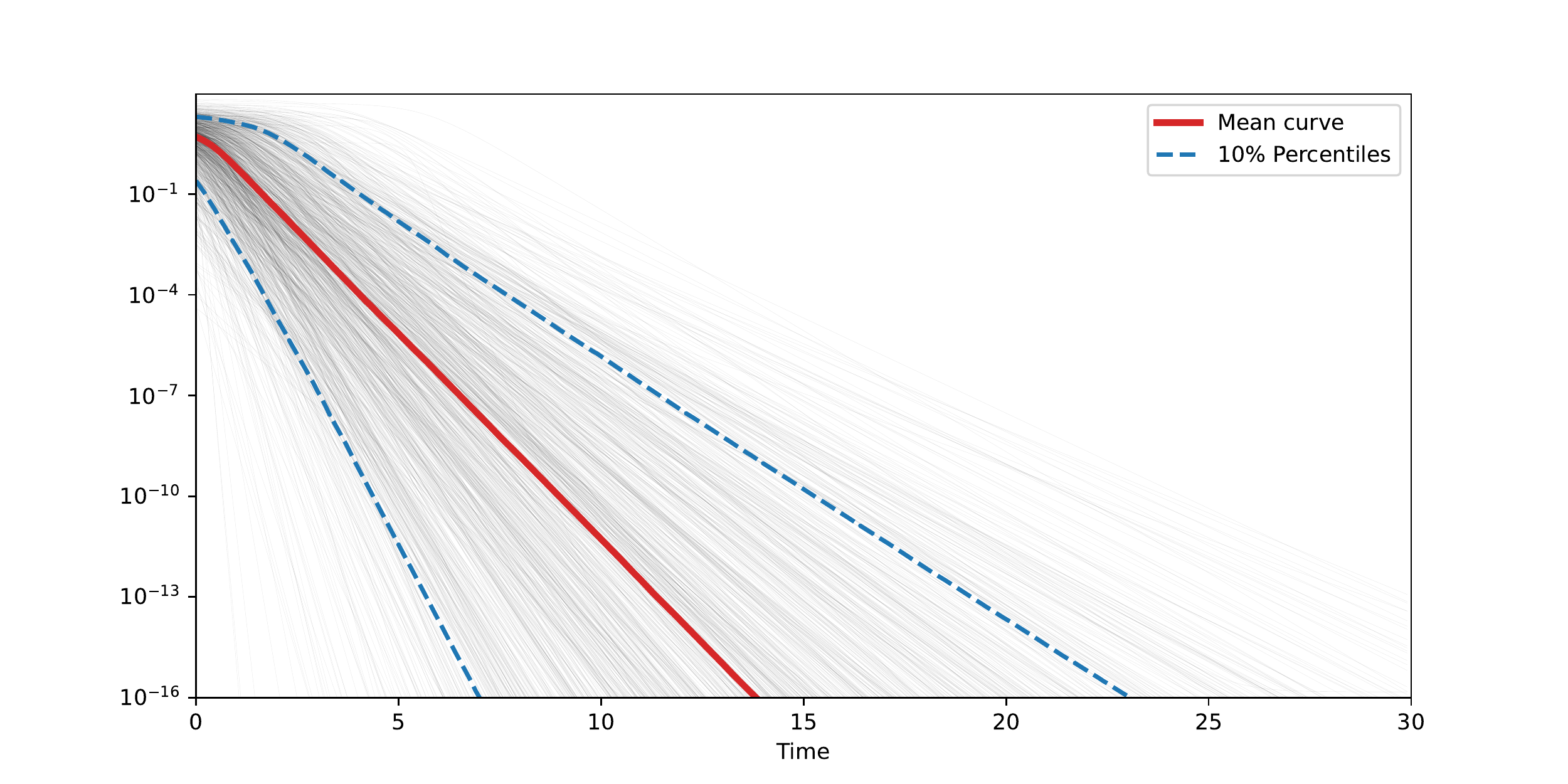}
    \caption{Plot illustrating the convergence of $1000$ gradient flow paths, each starting in a randomly selected matrix. The squared distance between each path and its corresponding polar decomposition is plotted (black curves) as well as the median squared distance (red curve) and percentiles (blue curves).}
    \label{fig:dist-conv}
\end{figure}

In Figure~\ref{fig:dist-conv}, the squared distance between each $B^{i}_k$ and $P^i$ is plotted against time, as well as the median, and upper and lower $10\%$ percentiles over all curves in each time step.
This illustrates the convergence of the gradient flow. 

Figure~\ref{fig:dist-conv} 
indicates exponential convergence. 
This hints at some intrinsic convexity of the problem. 

\section{Conclusion and outlook}

In the present paper, we have examined a gradient flow that computes the polar decomposition of matrices. 
It should be noted that if one just wants to compute the polar decomposition, the gradient flow is a poor choice. 
But, there are two main reasons for the study: First, as the introduction says, it shows how matrix decompositions and the Gaussian Monge problem have hidden dynamical features.

Second, here we have considered a finite-dimensional analogue of optimal transport in the smooth category, for which the same geometric structure holds, see \citet{Mo2017}. 
This paper may thus be seen as a pre-study for a detailed analysis of the vertical gradient flow in the infinite-dimensional setting, which could give rise to new algorithms for OMT.

 \bibliographystyle{amsplainnat}
 \bibliography{references}

\end{document}